\documentclass{article}

\usepackage{mathptmx}
\usepackage{eucal}
\usepackage{amsmath}
\usepackage{amscd}
\usepackage{amssymb}
\usepackage{amsthm}
\usepackage{xspace}
\usepackage[all,tips]{xy}
\usepackage[dvips]{graphicx}
\usepackage{verbatim}
\usepackage{syntonly}

\theoremstyle{plain}
\newtheorem{thm}{Theorem}[section]

\newtheorem{cor}[thm]{Corollary}

\newtheorem{prop}[thm]{Proposition}
\newtheorem{lemma}[thm]{Lemma}

\theoremstyle{definition}
\newtheorem*{rem}{Remark}
\newtheorem*{pf}{Proof}

\DeclareMathOperator{\Aut}{Aut} 
 \DeclareMathOperator{\Stab}{Stab}

\DeclareMathOperator{\rank}{rank}

 \DeclareMathOperator{\tr}{Tr}

\DeclareMathOperator{\SL}{SL} \DeclareMathOperator{\PSL}{PSL}
\DeclareMathOperator{\GL}{GL} 
\DeclareMathOperator{\Ort}{O}

 \DeclareMathOperator{\Sol}{Sol}
\DeclareMathOperator{\Nil}{Nil}

 \DeclareMathOperator{\Aff}{Aff}

\DeclareMathOperator*{\vol}{vol}

\DeclareMathOperator*{\Isom}{Isom}

\DeclareMathOperator*{\sign}{sign}

\newcommand{\eps}{\varepsilon}
\newcommand{\vp}{\varphi}

\newcommand{\al}{\alpha}

\newcommand{\be}{\beta}

\newcommand{\ga}{\gamma}
\newcommand{\Ga}{\Gamma}

\newcommand{\te}{\theta}

\newcommand{\si}{\sigma}

\newcommand{\de}{\delta}

\newcommand{\la}{\lambda}
\newcommand{\La}{\Lambda}

\newcommand{\ol}{\overline}

\newcommand{\fn}{\footnote}

\newcommand{\wt}{\widetilde}

\newcommand{\iny}{\infty}

\newcommand{\tri}{\ensuremath{\triangle}}
\newcommand{\prt}{\partial}

\newcommand{\co}{\ensuremath{\colon}}

\newcommand{\innp}[1]{\left< #1 \right>}
\newcommand{\abs}[1]{\left\vert#1\right\vert}
\newcommand{\set}[1]{\left\{#1\right\}}
\newcommand{\brac}[1]{\left[#1\right]}
\newcommand{\pr}[1]{\left( #1 \right) }

\newcommand{\su}{\subset}

\newcommand{\ba}{\bigcap}

\newcommand{\smin}{\setminus}

\newcommand{\bdef}{\overset{\text{def}}{=}}

\newcommand{\lra}{\longrightarrow}

\newcommand{\lmto}{\longmapsto}

\newcommand{\B}[1]{\ensuremath{\mathbf{#1}}}

\newcommand{\Cal}[1]{\ensuremath{\mathcal{#1}}}

\newcommand{\Hy}{\ensuremath{\B{H}}}
\newcommand{\N}{\ensuremath{\B{N}}}
\newcommand{\Q}{\ensuremath{\B{Q}}}
\newcommand{\R}{\ensuremath{\B{R}}}
\newcommand{\Z}{\ensuremath{\B{Z}}}
\newcommand{\C}{\ensuremath{\B{C}}}

\newcommand{\res}{\ensuremath{\textrm{Res}}}

\newcommand{\refP}[1]{Proposition~\ref{P:#1}}

\newcommand{\refSS}[1]{Subsection~\ref{SS:#1}}
\newcommand{\refT}[1]{Theorem~\ref{T:#1}}
\newcommand{\refL}[1]{Lemma~\ref{L:#1}}

\begin{document}


\title{\textbf{Cusps of Hilbert modular varieties}}
\author{D. B. McReynolds\fn{Supported by a V.I.G.R.E. graduate fellowship and Continuing Education fellowship.}}
\maketitle


\begin{abstract} \noindent Motivated by a question of Hirzebruch on the possible topological types of
cusp cross-sections of Hilbert modular varieties, we give a
necessary and sufficient condition for a manifold $M$ to be
diffeomorphic to a cusp cross-section of a Hilbert modular variety.
Specialized to Hilbert modular surfaces, this proves that every
$\Sol$ 3--manifold is diffeomorphic to a cusp cross-section of a
(generalized) Hilbert modular surface. We also deduce an obstruction
to  geometric bounding in this setting. Consequently, there exist
$\Sol$ 3--manifolds that cannot arise as a cusp cross-section of a
1--cusped nonsingular Hilbert modular surface.
\end{abstract}


\bibliographystyle{amsplain}


\section{Introduction}

\subsection*{Main results}

\noindent It is a classical problem in topology to decide whether or
not a closed $n$--manifold $M$ bounds. Hamrick and Royster
\cite{HamrickRoyster82} resolved this in the affirmative for flat
$n$--manifolds and Rohlin \cite{Rohlin51} for closed 3--manifolds.
However, beyond these two classes there are few other settings where
the story is nearly this complete. The introduction of geometry to a
topological problem provides additional structure which can lead to
new insight. This philosophy serves as motivation for the primary
concern of this article which is a geometric notion of bounding and
its specialization to infrasolv manifolds. \smallskip\smallskip

\noindent Let $k$ be a totally real number field with $[k:\Q]=n$,
$\Cal{O}_k$ the ring of integers of $k$, and $\si_1,\dots,\si_n$
denote the $n$ real embeddings of $k$. The group $\PSL(2;\Cal{O}_k)$
is an arithmetic subgroup of the $n$--fold product $(\PSL(2;\R))^n$
(see \cite{BorelChandra62}) via the embedding $\xi \lmto
(\si_1(\xi),\dots,\si_n(\xi))$ for $\xi \in \PSL(2;\Cal{O}_k)$.
Through this embedding, $\PSL(2;\Cal{O}_k)$ acts with finite volume
on the $n$--fold product of real hyperbolic planes $(\Hy^2)^n$. The
group $\PSL(2;\Cal{O}_k)$ is called \emph{the Hilbert modular
group}. More generally, we call any subgroup $\La$ of $\PSL(2;k)$
which is commensurable with $\PSL(2;\Cal{O}_k)$ a \emph{Hilbert
modular group} and the quotients $(\Hy^2)^n/\La$, \emph{Hilbert
modular varieties}. In the case that $k$ is a real quadratic number
field, these quotients are called \emph{Hilbert modular surfaces}.
For more on Hilbert modular surfaces, see \cite{Hirzebruch73} or
\cite{vanderGeer88}.\smallskip\smallskip

\noindent The primary focus of this article is cusp cross-sections
of Hilbert modular varieties. These infrasolv manifolds are virtual
$n$--torus bundles over $(n-1)$--tori where $[k:\Q]=n$ and $\rank
\Cal{O}_k^\times = n-1$. For brevity, we simply call these
\emph{virtual $(n,n-1)$--torus bundles}. Recall that an $n$--torus
bundle over an $m$--torus is the total space of a fiber bundle with
base manifold $T^m$ and fiber $T^n$. We call such manifolds simply
\emph{$(n,m)$--torus bundles}. We say that $N$ is a \emph{virtual
$(n,m)$--torus bundle} if $N$ is finitely covered by an
$(n,m)$--torus bundle.\smallskip\smallskip

\noindent In \cite{McReynolds04A}, cusp cross-sections of real,
complex, and quaternionic arithmetic hyperbolic $n$--orbifolds were
classified. In this article, we continue this theme by classifying
cusp cross-sections of Hilbert modular varieties. By taking the
quotient of the associated neutered space for the Hilbert modular
group $\La$, we obtain a compact Riemannian $2n$--orbifold whose
totally geodesic boundaries are the cusp cross-sections equipped
with metrics (defined up to scaling) coming from the associated
solvable Lie group.\smallskip\smallskip

\noindent Before stating our first classification result, we
introduce an additional piece of terminology.\smallskip\smallskip

\noindent For a totally real number field $k$, we say $\be \in k$ is
\emph{totally positive} if $\si_j(\be) >0$ for $j=1,\dots,n$. We
denote the set of totally positive elements and totally positive
integers by $k_+$ and $\Cal{O}_{k,+}$, and define the sets
$k_+^\times = k_+ \cap k^\times$, $\Cal{O}_{k,+}^\times =
\Cal{O}_k^\times \cap \Cal{O}_{k,+}$. We say that a virtual torus
bundle $N$ is \emph{$k$--defined} if there exists a faithful
representation $\rho\colon \pi_1(N) \lra k \rtimes k_+^\times$. If
in addition $\rho(\pi_1(N))$ is commensurable with $\Cal{O}_k
\rtimes \Cal{O}_{k,+}^\times$, we say that $N$ is
\emph{$k$--arithmetic}.\smallskip\smallskip

\noindent Our first result is:

\begin{thm}\label{T:Torus} A virtual $(n,n-1)$--torus bundle $N$ is diffeomorphic to a cusp
cross-section of a Hilbert modular variety over $k$ if and only if
$\pi_1(N)$ is $k$--arithmetic.
\end{thm}

\noindent \refT{Torus} answers a question of Hirzebruch \cite[page
203]{Hirzebruch73} who asked (in our terminology) which
$k$--arithmetic torus bundles arise as cusp cross-sections of
Hilbert modular varieties. See \refSS{Hirzebruch} for more on
this.\smallskip\smallskip

\noindent Every $(2,1)$--torus bundle admits either a Euclidean,
$\Nil$, or $\Sol$ geometry. Long and Reid \cite{LongReid02} proved
that the $(2,1)$--torus bundles which admit a Euclidean structure
are diffeomorphic to cusp cross-sections of arithmetic real
hyperbolic 4--orbifolds. In \cite{McReynolds04A}, we proved that
those that admit $\Nil$ structures are diffeomorphic to cusp
cross-sections of arithmetic complex hyperbolic 2--orbifolds. In
this article, we prove (see \S \ref{S:Sol} for the definitions):

\begin{thm}\label{T:21} Every $\Sol$ 3--manifold is diffeomorphic to a cusp cross-section
of a generalized Hilbert modular surface.
\end{thm}

\noindent We note that this shows closed 3--manifolds modelled on
this three geometries bound; of course, this is not new as Rohlin
proved this for any 3--manifold.\smallskip\smallskip

\noindent Using the Atiyah-Patodi-Singer signature formula, Long and
Reid \cite{LongReid00} showed that a flat 3--manifold which arises
as a cusp cross-section of a 1--cusped real hyperbolic 4--manifold
must have integral $\eta$--invariant. Together with Ouyang's work,
this proves that certain flat 3--manifolds cannot be the cusp
cross-section of a 1--cusped real hyperbolic 4--manifold. We
conclude this article with a similar result. Specifically, using the
work of Hirzebruch \cite{Hirzebruch73}, Atiyah-Donnely-Singer
\cite{ADS83}, and Cheeger-Gromov \cite{CheegerGromov85}, we prove:

\begin{thm}\label{T:GB} There exists a $\Sol$ 3--manifold which cannot be diffeomorphic
to a cusp cross-section of any 1--cusped Hilbert modular surface
with torsion free fundamental group.
\end{thm}

\subsection*{Acknowledgments}

\noindent I would like to thank my advisor Alan Reid for all his
help. In addition, I would like to thank Richard Schwartz for
suggesting Hilbert modular varieties as a family of examples for
which the techniques developed in \cite{McReynolds04A} might be
applied and for carefully reading an early draft of this article.

\section{Preliminary material}

\subsection{Stabilizer groups}

\noindent For $v \in \prt \Hy^n$, the group $\Stab(v) = \set{\ga \in
\Isom(\Hy^n)~:~\ga v = v}$ is isomorphic to $\R^{n-1} \rtimes (\R^+
\times \Ort(n-1))$. For $v \in \prt \Hy^n$ and $H<\Isom(\Hy^n)$, we
define the \emph{stabilizer group of $H$ at $v$} to be $\tri_v(H) =
\Stab(v) \cap H$. When $\tri_v(H)$ contains a parabolic isometry, we
call $\tri_v(H)$ the \emph{maximal peripheral subgroup of $H$ at
$v$} and say that $H$ has a \emph{cusp at $v$}. Often, we simply
write $\tri(H)$.\smallskip\smallskip

\noindent Cusps, horospheres, and cusp cross-sections are defined as
in the hyperbolic setting via Iwasawa decompositions of
$(\PSL(2;\R))^r$. For the Hilbert modular group $\PSL(2;\Cal{O}_k)$
over a totally real number field $k$, the stabilizer of the boundary
point corresponding to the Iwasawa decomposition given by the
$r$--fold product of the groups $\B{A},\B{N},\B{K}$ is the
peripheral subgroup
\[ \tri = \set{\begin{pmatrix} \be^{-1} & \al \\ 0 & \be \end{pmatrix}~:~\al
\in \Cal{O}_k,~\be \in \Cal{O}_{k,+}^\times}. \] Every peripheral
subgroup of $\PSL(2;\Cal{O}_k)$ is conjugate in $\PSL(2;k)$ to a
group commensurable with $\tri$.
\smallskip\smallskip

\subsection{Infrasolv manifolds and smooth rigidity}

\noindent For a simply connected, connected solvable Lie group $S$,
the affine group of $S$ is $\Aff(S) = S \rtimes \Aut(S)$. We say
that a discrete subgroup $\Ga < \Aff(S)$ is an \emph{infrasolv group
modelled on $S$} if $\Ga \cap S$ is finite index in $\Ga$ and
$S/\Ga$ is compact. An infrasolv group which is a subgroup of $S$
will be called a \emph{solv group modelled on $S$}. Any smooth
manifold which is diffeomorphic to $S/\Ga$ for some infrasolv group
will be called an \emph{infrasolv manifold modelled on $S$}. When
$\Ga$ is a solv group, we call the manifold $S/\Ga$ a \emph{solv
manifold modelled on $S$}.\smallskip\smallskip

\noindent We require the following rigidity result of Mostow
\cite{Mostow54}.

\begin{thm}[Mostow; \cite{Mostow54}]\label{T:Rig}  Let $M_1$ and $M_2$ be infrasolv manifolds.
If $\pi_1(M_1) \cong \pi_1(M_2)$, then $M_1$ is diffeomorphic to
$M_2$.
\end{thm}

\section{Cusps of Hilbert modular varieties}

\noindent In this section, we prove \refT{Torus}. The philosophy for
the proof is simple. Using the arithmeticity assumption on the torus
bundle $N$, we construct an injective homomorphism $\rho\co \pi_1(N)
\lra \tri(\PSL(2;\Cal{O}_k))$. To find a Hilbert modular group $\La$
for which $\tri(\La) = \rho(\pi_1(N))$, we are reduced to making a
subgroup separability argument. The proof is completed by applying
\refT{Rig}. The remainder of this section is devoted to the details.

\subsection{Subgroup separability}

\noindent Recall that if $G$ is a group, $H<G$ and $g \in G \smin
H$, we say $H$ and $g$ are \emph{separated} if there exists a
subgroup $K$ of finite index in $G$ which contains $H$ but not $g$.
We say that $H<G$ is \emph{separable} in $G$ if every $g \in G \smin
H$ and $H$ can be separated.\smallskip\smallskip

\noindent As in \cite{McReynolds04A}, the main technical result we
make use of is:

\begin{thm}\label{T:Borel} Let $\La$ be a Hilbert modular group and $\tri(\La)$,
a maximal peripheral subgroup. Then every subgroup of $\tri(\La)$ is
separable in $\La$.
\end{thm}

\subsection{Proof of \refT{Torus}}

\noindent In this subsection, we prove \refT{Torus}. The following
establishes a correspondence between $k$--arithmetic torus bundle
groups and maximal peripheral subgroups of Hilbert modular groups.

\begin{thm}[Correspondence theorem]\label{T:CT} Let $N$ be a $k$--arithmetic torus bundle.
Then there exists a faithful representation $\psi\co \pi_1(N) \lra
\tri(\PSL(2;\Cal{O}_k))$ such that $\psi(\pi_1(N))$ is a finite
index subgroup of $\tri(\PSL(2;\Cal{O}_k))$. Moreover, there exists
a finite index subgroup $\La$ of $\PSL(2;\Cal{O}_k)$ such that
$\tri(\La) = \psi(\pi_1(N))$.
\end{thm}

\noindent We defer the proof of \refT{CT} for the moment in order to
prove \refT{Torus}. \smallskip\smallskip

\begin{pf}[Proof of \refT{Torus}] For the direct implication, since $N$ is diffeomorphic to a cusp
cross-section of a Hilbert modular variety, there exists a Hilbert
modular group $\La$ and an isomorphism $\psi\co \pi_1(N) \lra
\tri(\La)$. To obtain an injective homomorphism $\rho\co \pi_1(N)
\lra k \rtimes k_+^\times$ such that $\rho(\pi_1(N))$ is
commensurable with $\Cal{O}_k \rtimes \Cal{O}_{k,+}^\times$, we
argue as follows. By conjugating by an element $\ga$ of $\PSL(2;k)$,
we can assume that
\[ \ga^{-1}\psi(\pi_1(N))\ga \su B_k = \set{\begin{pmatrix} \be^{-1} &
\al \\ 0 & \be \end{pmatrix} ~:~ \al \in k,~\be \in k_+^\times}. \]
As $\ga \in \PSL(2;k)$, $\ga^{-1}\La\ga$ remains a Hilbert modular
group, and moreover, $\ga^{-1}\psi(\pi_1(N))\ga$ is commensurable
with
\[ \tri(\PSL(2;\Cal{O}_k)) = \set{\begin{pmatrix} \be^{-1} & \al \\ 0
& \be \end{pmatrix} ~:~ \al \in \Cal{O}_k,~\be \in
\Cal{O}_{k,+}^\times}. \] To obtain the faithful representation
$\rho$, we simply compose $\mu_\ga \circ \psi$ with the isomorphism
$\iota\co B_k \lra k \rtimes k_+^\times$ given by
$\iota\pr{\pr{\begin{smallmatrix} \be^{-1} & \al \\ 0 & \be
\end{smallmatrix}}} = (\al,\be)$.\smallskip\smallskip

\noindent For the reverse implication, we apply \refT{CT} and
\refT{Rig}. Specifically, let $\La$ by the Hilbert modular group
guaranteed by \refT{CT} and let $N^\prime$ denote an embedded cusp
cross-section associated with $\tri(\La)$. As a smooth manifold,
$N^\prime$ is of the form $\R^{2n-1}/\tri(\La)$. By \refT{CT}, we
have an isomorphism $\psi\co \pi_1(N) \lra \pi_1(N^\prime)$.
Applying \refT{Rig}, we obtain the desired diffeomorphism between
$N$ and $N^\prime$.
\end{pf}

\noindent In the proof of \refT{CT}, the following lemma is
required.

\begin{lemma}\label{L:1} Let $N$ be a $k$--arithmetic torus bundle.
Then there exists an injective homomorphism $\rho\colon \pi_1(N)
\lra \Cal{O}_k \rtimes \Cal{O}_{k,+}^\times$. Moreover,
$\rho(\pi_1(N))$ is a finite index subgroups of $\Cal{O}_k \rtimes
\Cal{O}_{k,+}^\times$.
\end{lemma}

\begin{pf} Since $N$ is $k$--arithmetic, we have a faithful representation
$\te\colon \pi_1(N) \lra k \rtimes k_+^\times$ such that
$\te(\pi_1(N))$ is commensurable with $\Cal{O}_k \rtimes
\Cal{O}_{k,+}^\times$. Hence, given $(\al,\be)\in \te(\pi_1(N))$, we
have for some $m \in \N$,
\[ (\al + \be\al + \be^2\al + \dots+ \be^{m-1}\al,\be^m) \in \Cal{O}_k
\rtimes \Cal{O}_{k,+}^\times. \] Consequently, $\be^m \in
\Cal{O}_{k,+}^\times$ and thus $\be \in \Cal{O}_{k,+}^\times$. Even
so, it may be the case that $(\al,\be)$ is not contained in
$\Cal{O}_k \rtimes \Cal{O}_{k,+}^\times$. This is rectified as
follows. Select a generating set for $\pi_1(N)$, say
$g_1,\dots,g_u$. For each generator, we have $\te(g_j) =
(\al_j,\be_j)$ with $\al_j\in k$ and $\be_j \in
\Cal{O}_{k,+}^\times$. Since $k$ is the field of fractions of
$\Cal{O}_k$,  we can select $\la_j \in \Cal{O}_k$ such that
$(0,\la_j)\te(g_j)(0,\la_j)^{-1} \in \Cal{O}_k \rtimes
\Cal{O}_{k,+}^\times$. Note that
\[ (0,\la_j)\te(g_j)(0,\la_j)^{-1} = (\la_j\al_j,\be_j), \]
and so the second coordinate $\be_j$ is unchanged. Finally, for $\la
= \la_1\dots\la_u$, define $\rho = \mu_{(0,\la)} \circ \te$, where
$\mu_{(0,\la)}$ denotes the inner automorphism determined by
$(0,\la)$. By construction, $\rho$ is a faithful representation of
$\pi_1(N)$ onto a finite index subgroup of $\Cal{O}_k \rtimes
\Cal{O}_{k,+}^\times$.
\end{pf}

\noindent With \refL{1} in hand, we prove
\refT{CT}.\smallskip\smallskip

\begin{pf}[Proof of \refT{CT}] By \refL{1}, we have an injective homomorphism $\rho\co
\pi_1(N) \lra \Cal{O}_k \rtimes \Cal{O}_{k,+}^\times$ such that
$\rho(\pi_1(N))$ is a finite index subgroup. To obtain the injective
homomorphism $\psi$, we compose $\rho$ with the isomorphism
\[ \iota^{-1}\co \Cal{O}_k \rtimes \Cal{O}_{k,+}^\times \lra \tri(\PSL(2;\Cal{O}_k)) \]
where $\iota^{-1}(\al,\be) = \pr{\begin{smallmatrix} \be^{-1} & \al
\\ 0 & \be
\end{smallmatrix}}$. That $\psi$ is faithful and $\psi(\pi_1(N))$ is a finite index
subgroup of $\tri(\PSL(2;\Cal{O}_k))$ follow immediately from the
properties of $\rho$ and $\iota$. \smallskip\smallskip

\noindent To find the desired subgroup $\La$, we apply \refT{Borel}.
Specifically, select a complete set of coset representatives
$\ga_1,\dots,\ga_s$ for $\tri(\PSL(2;\Cal{O}_k)) /\psi(\pi_1(N))$.
By \refT{Borel}, $\psi(\pi_1(N))$ is separable. Therefore for each
$j$ we can find finite index subgroups $\La_j$ such that $\ga_j
\notin \La_j$ and $\psi(\pi_1(N)) < \La_j$. To get the desired
$\La$, take $\La = \ba_{j=1}^s \La_j$.
\end{pf}

\subsection{A question of Hirzebruch}\label{SS:Hirzebruch}

\noindent Let $k$ be a totally real number field, $M<k$ an additive
group of rank $n$ (the degree of $k$ over $\Q$), and
$V<\Cal{O}_{k,+}^\times$ a finite index subgroup such that for all
$\la \in V$, $\la M \su M$. For each pair $(M,V)$, we define the
peripheral group
\[ \tri(M,V) = \set{\begin{pmatrix} \be^{-1} & \al \\ 0 & \be
  \end{pmatrix} ~:~ \al \in M,~\be \in V} < \PSL(2;k). \]
For any Hilbert modular variety, the peripheral groups $\tri(\La)$
are conjugate (in $\PSL(2;k)$) to groups of the form $\tri(M,V)$. In
\cite[p. 203]{Hirzebruch73}, Hirzebruch mentions that it is
apparently unknown whether or not every $\tri(M,V)$ can occur as a
maximal peripheral subgroup of a Hilbert modular group. The
following corollary gives an affirmative answer.

\begin{cor}\label{C:Hirzebruch} For every pair $(M,V)$, there exists a
Hilbert modular group $\La$ such that $\tri(\La) = \tri(M,V)$.
\end{cor}

\begin{pf} As in the proof of \refL{1}, we can conjugate $\tri(M,V)$ by an
element of the form $\ga = \pr{\begin{smallmatrix} \la^{-1} & 0 \\ 0
& \la \end{smallmatrix}}$, with $\la \in \Cal{O}_k$, such that
$\ga^{-1} \tri(M,V) \ga$ is contained in $\PSL(2;\Cal{O}_k)$. Since
$M$ and $V$ are finite index subgroups of $\Cal{O}_k$ and
$\Cal{O}_{k,+}^\times$, respectively, $\ga^{-1}\tri(M,V)\ga$ is a
finite index subgroup of $\tri(\PSL(2;\Cal{O}_k))$. Thus there
exists a finite index subgroup $\La_1<\PSL(2;\Cal{O}_k)$ such that
$\tri(\La_1) = \ga^{-1}\tri(M,V)\ga$. Hence, for $\La = \ga \La_1
\ga^{-1}$, we have $\tri(\La) = \tri(M,V)$. As $\ga \in \PSL(2;k)$,
$\La$ is a Hilbert modular group, as required.
\end{pf}

\section{A simple criterion for arithmeticity}

\noindent In this section, we give a simple criterion for the
arithmeticity of $(n,m)$--torus bundles. The need for such a result
is practical, as it allows one to establish the arithmeticity of a
torus bundle computationally. We encourage the reader to compare the
results of this section with Corollary 5.5 in \cite{McReynolds04A}.

\subsection{Linear equations and presentations of torus bundle groups}

\noindent For an (orientable) $(n,n-1)$--torus bundle $M$, since
both the base and fiber are aspherical, we have the short exact
sequence induced by the long exact sequence of the fiber bundle
\[ 1 \lra \Z^n \lra \pi_1(M) \lra \Z^{n-1} \lra 1. \]
The action of $\Z^{n-1}$ on $\Z^n$ induces a homomorphism $\vp\colon
\Z^{n-1} \lra \SL(n;\Z)$ called the \emph{holonomy representation}.
Since peripheral subgroups in Hilbert modular groups have faithful
holonomy representation, we assume throughout that $\vp$ is
faithful. In particular, we obtain a faithful representation of
$\pi_1(M)$ into $\Z^n \rtimes \SL(n;\Z)$.\smallskip\smallskip

\noindent Of primary importance for us here is that the holonomy
representation together with any finite presentation yields a
homogenous linear system of equations with coefficients in $\Z$.
This system arises as follows. For ease, select a presentation of
the form
\[ \innp{x_1,\dots,x_n,\ol{y_1},\dots,\ol{y_{n-1}}~:~R} \]
where $x_1,\dots,x_n$ generate $\Z^m$, $\ol{y_1},\dots,\ol{y_{n-1}}$
are lifts of a generating set $y_1,\dots,y_{n-1}$ for $\Z^{n-1}$,
and $R$ is a finite set of relations of the form
\[ x_j\ol{y_k} = \ol{y_k}w_{j,k}, \quad w_{j,k} \in
\innp{x_1,\dots,x_n}. \] Using the holonomy representation, we can
write
\[ x_j = (a_j,I), \quad \ol{y_j} = (b_j,\vp(y_j)) \in \R^n \rtimes
\SL(n;\R). \] Each relation in the presentation yields a linear
homogenous equation in the vector variables $a_j$ and $b_j$ (see
below for an explicit example of how these equations arise). Namely,
we insert the above forms for $x_j$ and $\ol{y_k}$ into the relation
and consider only the first coordinate. The equations we obtain are
of the form
\[ a_j + b_k - \vp(y_k) - v_{j,k} = 0 \]
where $w_{j,k} = (v_{j,k},I)$. That this system has integral
solutions which yield faithful representations follows from the fact
that $\vp$ is faithful and induces a faithful representation of
$\pi_1(M)$ into $\Z^n \rtimes \SL(n;\Z)$.

\subsection{A simple criterion for arithmeticity}

\noindent The main result of this section is a simple criterion for
arithmeticity based on the structure of the holonomy representation.
In the statement and proof, let $\res_{k/\Q}$ denotes restriction of
scalars from $k$ to $\Q$ and assume that $[k:\Q]=n$ and $\rank
\Cal{O}_k^\times = n-1$. In particular, $k$ is totally real.

\begin{thm}\label{T:Class} Let $M$ be an orientable $(n,n-1)$--torus bundle. Then $M$ is diffeomorphic to a cusp
cross-section of a Hilbert modular variety defined over $k$ if and
only if $\vp = \res_{k/\Q}(\chi)$, for some faithful character
$\chi\colon \Z^{n-1} \lra \Cal{O}_{k,+}^\times$, where $\vp$ is some
holonomy representation.
\end{thm}

\begin{pf} For the direct implication, since $M$ is diffeomorphic to a cusp
cross-section of a Hilbert modular variety, by \refT{Torus}, we have
a faithful representation
\[ \rho\co \pi_1(M) \lra \Cal{O}_k \rtimes \Cal{O}_{k,+}^\times \]
By restricting scalars from $k$ to $\Q$, we obtain a faithful
representation
\[ \res_{k/\Q}(\rho)\co \pi_1(M) \lra \Z^n \rtimes \SL(n;\Z). \]
The proof is completed by noting that the holonomy map induced by
this representation is simply $\res_{k/\Q}(\chi)$, where $\chi\co
\Z^{n-1} \lra \Cal{O}_{k,+}^\times$ is the holonomy representation
induced by the representation $\rho$.\smallskip\smallskip

\noindent For the converse, we seek a faithful representation
$\rho\colon \pi_1(M) \lra \Cal{O}_k \rtimes \Cal{O}_{k,+}^\times$.
Note that since $[k:\Q]=n$ and $\rank \Cal{O}_k^\times = n-1$, the
image of $\pi_1(M)$ would necessarily be a finite index subgroup. By
assumption, we have a faithful character $\chi\colon \Z^{n-1} \lra
\Cal{O}_{k,+}^\times$. We extend this to a faithful representation
of $\pi_1(M)$ into $\Cal{O}_k \rtimes \Cal{O}_{k,+}^\times$ as
follows. Select a presentation as above for $\pi_1(M)$ with
generators $x_1,\dots,x_n,\ol{y_1},\dots, \ol{y_{n-1}}$. Write
\begin{equation}\label{E:Hom}
x_i = (\al_i,1),~\ol{y_i} = (\ga_i,\chi(y_i)) \in k \rtimes
\Cal{O}_{k,+}^\times
\end{equation}
where $\al_i$ and $\ga_i$ are to be determined. Using our
presentation for $\pi_1(M)$, we obtain a system of linear homogenous
equations $\Cal{L}$ with coefficients in $\Cal{O}_k$. Note, as
above, solutions to $\Cal{L}$ yield representations of $\pi_1(M)$
into $k\rtimes \Cal{O}_{k,+}^\times$. We assert that there is a
solution which yields a faithful representation. To see this, by
restricting scalars from $k$ to $\Q$, we obtain a linear system
$\res_{k/\Q}(\Cal{L})$ with coefficients in $\Z$. Solutions to the
system $\res_{k/\Q}(\Cal{L})$ yield representations of $\pi_1(M)$
into $\Z^n \rtimes \SL(n;\Z)$. Moreover, a solution to
$\res_{k/\Q}(\Cal{L})$ which yields a faithful representation is
equivalent to a solution of $\Cal{L}$ which yields a faithful
representation into $\Cal{O}_k \rtimes \Cal{O}_{k,+}^\times$. That
such a solution exists with integral coefficients for
$\res_{k/\Q}(\Cal{L})$ follows from the faithfulness of
$\res_{k/\Q}(\chi)$ and our discussion in the previous subsection.
This yields a solution for $\Cal{L}$ with coefficients in
$\Cal{O}_k$ which yields a faithful representation. Therefore, $M$
is $k$--arithmetic, since there exists a faithful representation
$\psi\co \pi_1(M) \lra \Cal{O}_k \rtimes \Cal{O}_{k,+}^\times$ such
that $\psi(\pi_1(M))$ is a finite index subgroup of $\Cal{O}_k
\rtimes \Cal{O}_k^\times$.
\end{pf}

\begin{rem} If the character $\chi$ only maps into $\Cal{O}_k^\times$, the above
proof yields a faithful representation $\rho\colon \pi_1(M) \lra
\Cal{O}_k \rtimes \Cal{O}_k^\times$.
\end{rem}

\section{$\Sol$ 3--manifolds}\label{S:Sol}

\noindent Before proving \refT{21}, we give a brief review of $\Sol$
3--manifolds (see \cite{Scott83}). Let $\Sol= \R^2 \times \R^+$ with
group operation defined by
\[ (x_1,y_1,t_1)\cdot(x_2,y_2,t_2) \bdef (x_1 + e^{t_1}x_2,y_1+
e^{-t_1}y_2, t_1+t_2). \] By a \emph{$\Sol$ 3--orbifold}, we mean a
manifold $M$ which is diffeomorphic to $\Sol/\Ga$, where $\Ga$ is a
discrete subgroup of $\Aff(\Sol)$ such that $\Sol/\Ga$ is compact
and $[\Ga:\Ga \cap \Sol]<\iny$. These manifolds, in the terminology
from \S 2, are infrasolv manifolds modelled on $\Sol$. However, the
terminology used in this section for these manifolds is more
prevalent. \smallskip\smallskip

\noindent In \cite{Scott83}, Scott proved that every $(2,1)$--torus
bundles admits either a Euclidean, $\Nil$, or $\Sol$ structure. The
following result is easily derived from  \cite{Scott83}. We include
a proof here for completeness.

\begin{prop}\label{P:Scott} Let $M$ be an orientable $(2,1)$--torus bundle
which admits a $\Sol$ structure. Then there exists a faithful
representation $\rho\colon \pi_1(M) \lra \Cal{O}_k \rtimes
\Cal{O}_k^\times$ for some real quadratic number field $k$.
\end{prop}

\begin{pf} For any $(2,1)$--torus bundle $M$, let the $\Z$--action be given by
$A = \pr{\begin{smallmatrix} a & b \\ c & d \end{smallmatrix}}$. If
the order of $A$ is finite, then $\pi_1(M)$ is a Bieberbach group
and $M$ admits a Euclidean structure. Therefore we may assume that
the order of $A$ is infinite. If $A$ is not diagonalizable, then
some power of $A$ is conjugate to $\pr{\begin{smallmatrix} 1 & \al
\\ 0 & 1 \end{smallmatrix}}$ with $\al \ne 0$. In this case, $M$
admits a $\Nil$ structure. Thus, we may assume that $A$ is
diagonalizable. In this case we have $\pr{\begin{smallmatrix} \be &
0 \\ 0 & \be^{-1} \end{smallmatrix}}$ for a conjugate of $A$. It
follows, since $A \in \SL(2;\Z)$, that $\be$ and $\be^{-1}$ are
algebraic integers in the real quadratic field $\Q(\be)$. Thus the
representation $\vp\colon \Z \lra \GL(2;\Z)$ is conjugate to
$\res_{k/\Q}(\chi)$, where $\chi\colon \Z \lra \Cal{O}_k^\times$ is
given by $\chi(1) = \be$. Therefore by the remark following
\refT{Class}, we have a faithful representation $\rho\colon \pi_1(M)
\lra \Cal{O}_k \rtimes \Cal{O}_k^\times$, as asserted.
\end{pf}

\noindent Via \refP{Scott}, note every $\Sol$ 3--manifold group does
faithfully represent into $\Isom((\Hy^2)^2)$. Those that arise as
cusp cross-sections of Hilbert modular surfaces are precisely the
ones whose fundamental group faithfully represents into the identity
component of $\Isom((\Hy^2)^2)$. However, the quotients of those
groups which fail to map into the identity component do produce
finite volume quotients which possess 2--fold covers which are
Hilbert modular surfaces. For this reason, we call such quotients
\emph{generalized Hilbert modular varieties}. Given this, \refT{21}
follows from this discussion in combination with \refT{Torus}.

\section{Geometric bounding}

\noindent Let $W$ be a 1--cusped Hilbert modular surface $W$ with
torsion free fundamental group---we call $W$ a \emph{Hilbert modular
manifold} in this case. Similar to the thick-thin decomposition of a
real hyperbolic $n$--manifold, $W$ has a decomposition comprised of
a compact manifold $\wt{W}$ with boundary $S$ and cusp end $S \times
\R^+$. Following Schwartz \cite{Schwartz95} (see also
\cite{FarbSchwartz96}), we call the universal cover of $\wt{W}$ the
associated \emph{neutered manifold}, and note $\wt{W}$ is a compact
4--manifold with $\Sol$ 3--manifold boundary. Moreover, the locally
symmetric metric $\wt{g}$ on $W$ restricted to $S$ endows $S$ with a
complete $\Sol$ metric $g$ such that $\wt{g}$ is a complete, finite
volume metric in the interior of $\wt{W}$ and $(S,g)$ is a totally
geodesic boundary.\smallskip\smallskip

\noindent The goal of this section is the establishment of a
nontrivial obstruction for this geometric situation. The obstruction
is obtained by mimicking the argument of Long--Reid
\cite{LongReid00} for flat 3--manifolds. This in combination with a
calculation of Hirzebruch bears \refT{GB} from the
introduction.\smallskip\smallskip

\noindent In \cite{Hirzebruch73}, Hirzebruch extended his signature
formula to Hilbert modular surfaces. The formula relates the
signature of the neutered manifold $\wt{W}$ to a Hirzebruch
$L$--polynomial evaluated on the Pontrjagin classes of $\wt{W}$ but
with a correction term associated to $\prt \wt{W}$. When $\pi_1(W)$
contains torsion, the elliptic singularities also contribute
nontrivially to this correction term, and so for simplicity, we
assume throughout that $\pi_1(W)$ is torsion free. In this case,
Hirzebruch's formula becomes
\[ \si(\wt{W}) = \de(E_1)+\dots +\de(E_r) \]
where $E_1,\dots,E_r$ is a complete set of cusp ends of $W$ given
from the thick-thin decomposition and $\si(\wt{W})$ denotes the
signature of $\wt{W}$. The definition of the terms $\de(E_j)$ are
given as follows. Associated to each cusp end is the
$\pi_1(W)$--conjugacy class of a maximal peripheral subgroup
$\Ga_j$. The group $\Ga_j$ is conjugate in $\PSL(2;k)$ to a subgroup
of the familiar form $\tri(M_j,V_j)$. In turn, for the pair
$(M_j,V_j)$, we have an associated Shimuzu $L$--function
$L(M_j,V_j,s)$---see \cite{Shimizu63}---defined by
\[ L(M,V,s) = \sum_{\be \in (M_j\smin \set{0})/V_j} \frac{\sign(N_{k/\Q}(\be))}{(N_{k/\Q}(\be))^s} \]
where $N_{k/\Q}$ is the norm map. With this, the invariant
$\de(E_j)$ is defined to be
\[ \de(E_j) = \frac{-\vol(M_j)}{\pi^2}L(M_j,V_j,1) \]
where $\vol(M_j)$ is the volume of $\R^2/M$ with respect to the
pairing $\tr_{k/\Q}$. Equivalently,
\[ \vol(M_j) = \abs{\det(\be_i^{(j)})}, \]
where $\be_1,\be_2$ is a $\Z$--module basis for $M_j$ and
$\be_i^{(1)}$ and $\be_i^{(2)}$ denote the image of $\be_i$ under
the two real embeddings of $k$ into $\R$.

\begin{thm}[Hirzebruch;\cite{Hirzebruch73}]\label{T:Signature}
If $W$ is a Hilbert modular manifold with exactly one cusp, then
\[ \si(\wt{W}) = \frac{-\vol(M)}{\pi^2}L(M,V,1) \]
for the unique $\pi_1(W)$--conjugacy class $\tri(M,V)$.
\end{thm}

\noindent As we seek an integrality condition, it is convenient to
change the pair $M,V$. Associated to the $\Z$--module $M$ is the
\emph{dual lattice} $M^*$ defined to be the image of $M$ under the
duality pairing provided by $\tr_{k/\Q}$.

\begin{prop}\label{P:Dual}
For a horosphere $\Cal{H}$ stabilized by $\tri(M,V)$ and
$\tri(M^*,V)$, $\Cal{H}/\tri(M,V)$ and $\Cal{H}/\tri(M^*,V)$ are
diffeomorphic $\Sol$ 3--manifolds.
\end{prop}

\begin{pf} Let $\vp_M,\vp_{M^*}\co V \lra \SL(2;\Z)$
be the holonomy representations for $\tri(M,V)$ and $\tri(M^*,V)$.
The pairing $\tr_{k/\Q}$ can be viewed as an element of $\la
\in\SL(2;\Z)$ such that $\la M = M^*$. By construction $\vp_{M^*} =
\la(\vp_M)\la^{-1}$, and so we have an isomorphism $\rho\co
\tri(M,V) \lra \tri(M^*,V)$ given by
\[ \rho(\be,\vp_M(\al)) = (\la \be, \la\vp_M(\al)\la^{-1}). \]
The proof is completed by appealing to the smooth rigidity theorem
of Mostow \refT{Rig}.
\end{pf}

\noindent Hecke (see \cite{ADS83}) related the $L$--functions
$L(M,V,s)$ and $L(M^*,V,s)$ by the functional equation $H(M,V,s) =
(-1)^sH(M^*,V,1-s)$, where
\[ H(M,V,s) = \brac{\Ga\pr{\frac{s+1}{2}}}^2\pi^{-(s+1)}\brac{\vol(M)}^sL(M,V,s) \]
The specialization of this functional equation at $s=1$ produces
\begin{align*}
\pr{\Ga(1)}^2\pi^{-2}\vol(M)L(M,V,1) &= -\pr{\Ga\pr{\frac{1}{2}}}^2\pi^{-1}L(M^*,V,0) \\
L(M^*,V,0) &= -\frac{\vol(M)}{\pi^2}L(M,V,1),
\end{align*}
and thus from this and \refT{Signature}, we obtain
\begin{equation}\label{e:eta}
\si(\wt{W}) = L(M^*,V,0).
\end{equation}

\noindent It is at this point that we take stock in what has been
done. For a 1--cusped Hilbert modular manifold $W$ with cusp
cross-section $S$, we have associated to $S$ the invariant $\de(S
\times \R^+)$. As both $M$ and $V$ depend on the associated $\Sol$
metric on $S$ afforded by its embedding as a cusp cross-section, the
invariant $\de(S\times \R^+)$ depends on the associated $\Sol$
metric on $S$. Our goal is to use the integrality of $\si(\wt{W})$
and (\ref{e:eta}) to produce an obstruction for $S$ to topologically
occur in this geometric setting. For this, it remains to show the
invariant $\de(S \times \R^+)$ is independent of the $\Sol$
structure on $S$.  \smallskip\smallskip

\noindent Given a peripheral group $\tri(M,V)$ and stabilized
horosphere $\Cal{H}$, the metric on $\Hy_\R^2 \times \Hy_\R^2$
endows $\Cal{H}$ with a $\tri(M,V)$--invariant metric
$g_{\Cal{H},M,V}$. Consequently the metric $g_{\Cal{H},M,V}$
descends to quotient $\Cal{H}/\tri(M,V)$ and endows
$\Cal{H}/\tri(M,V)$ with a complete $\Sol$ structure that depends on
the horosphere $\Cal{H}$ only up to similarity. \smallskip\smallskip

\noindent The formula (\ref{e:eta}) was also established in
\cite{ADS83} where $L(M^*,V,0)$ was reinterpreted as the
$\eta$--invariant of an adiabatic limit.

\begin{thm}[Atiyah--Donnely--Singer;\cite{ADS83}]\label{T:ADS1}
\[ L(M^*,V,0) = \lim_{\eps \lra 0} \eta(\Cal{H}/\tri(M^*,V),g_{\Cal{H},M^*,V}/\eps). \]
\end{thm}

\noindent More generally, given any $\Sol$ structure $g$ on $S$, we
can define
\[ \de(S,g) = \lim_{\eps \lra 0} \eta(S,g/\eps). \]
The last ingredient for proof of \refT{GB} is the independence of
$\de(S,g)$ from $g$, a result established by Cheeger and Gromov
\cite{CheegerGromov85}.

\begin{thm}[Cheeger--Gromov;\cite{CheegerGromov85}]\label{T:Eta}
$\de(S,g)$ is a topological invariant of the $\Sol$ 3--manifold $S$.
\end{thm}

\noindent We are now in position to state and prove the principal
observation needed in the proof of \refT{GB} (compare with
\cite{LongReid00}).

\begin{thm}\label{T:GeometricBounding}
If $S$ is diffeomorphic to a cusp cross-section of a 1--cusped
Hilbert modular manifold, then $\de(S) \in \Z$.
\end{thm}

\begin{pf}
If $(S,g)$ arises as a cusp cross-section of a 1--cusped Hilbert
modular manifold $W$, then there is an isometric embedding $f\co
(S,g) \lra W$ onto a cusp cross-section of $W$. Let $f_*(\pi_1(S))
=\tri(M,V)$ with associated horosphere $\Cal{H}$ selected such that
$\Cal{H}/\tri(M,V)$ is embedded in $W$. By \refP{Dual},
$\Cal{H}/\tri(M^*,V)$ is diffeomorphic to $S$, though equipped with
the metric $g_{\Cal{H},M^*,V}$. From the computation above in
combination with \refT{ADS1}, $\si(\wt{W}) =
\de(S,g_{\Cal{H},M^*,V})$ and by \refT{Eta}, the right hand side
depends only on the topological type of $S$. Since $\si(\wt{W})$ is
in $\Z$, $\de(S)$ is in $\Z$ as asserted.
\end{pf}

\begin{pf}[Proof of \refT{GB}]
To prove \refT{GB}, by \refT{GeometricBounding}, it suffices to find
a $\Sol$ 3--manifold $S$ for which $\de(S) \notin \Z$. For
$k=\Q(\sqrt{3})$, the standard Hilbert modular surface $W$ over $k$
has precisely one cusp, since the number of cusps of a standard
Hilbert modular surface over $k$ is the ideal class number of $k$.
Setting $S$ to be an embedding cusp cross-section of $W$, the proof
is completed by appealing to \cite{Hirzebruch73}. Specifically,
Hirzebruch showed $\de(S) = -1/3$.
\end{pf}

\begin{rem} It is unknown to the author whether or not there exist 1--cusped
Hilbert modular manifolds. In addition, the number fields
$\Q(\sqrt{6})$, $\Q(\sqrt{21})$ and $\Q(\sqrt{33})$ also have
standard Hilbert modular surfaces with precisely one cusp for which
the associated invariant $\de(S) \notin \Z$. In each of these cases,
$\de(S)=-2/3$ (see \cite[p. 236]{Hirzebruch73}).
\end{rem}

\noindent It is unknown to the author whether or not there exist
1--cusped Hilbert modular manifolds. Using the generalized Riemann
hypothesis, K. Petersen \cite{Petersen05} constructed infinite many
1--cusped Hilbert modular surfaces. However, the nature of the
construction likely produces Hilbert modular surface groups with
2--torsion.


\def\cprime{$'$} \def\lfhook#1{\setbox0=\hbox{#1}{\ooalign{\hidewidth
  \lower1.5ex\hbox{'}\hidewidth\crcr\unhbox0}}} \def\cprime{$'$}
  \def\cprime{$'$}
\providecommand{\bysame}{\leavevmode\hbox
to3em{\hrulefill}\thinspace}
\providecommand{\MR}{\relax\ifhmode\unskip\space\fi MR }
\providecommand{\MRhref}[2]{%
  \href{http://www.ams.org/mathscinet-getitem?mr=#1}{#2}
} \providecommand{\href}[2]{#2}


\noindent
Department of Mathematics, \\
California Institute of Technology, \\
253-37 Mathematics, Pasadena, CA 91125\\
email: {\tt dmcreyn@caltech.edu}


\end{document}